\documentclass[11pt,twoside]{article}
\usepackage{amsfonts,amsmath}

\newtheorem{definition}{Definition}[section]
\newenvironment{defi}{\begin{definition} \rm}{\end{definition}}

\newtheorem{prop}[definition]{Proposition}

\newtheorem{lemm}[definition]{Lemma}
\newtheorem{fact}[definition]{Fact}

\newtheorem{coro}[definition]{Corollary}
\newtheorem{theo}[definition]{Theorem}
\newtheorem{notation}[definition]{Notation}

\newtheorem{construction}[definition]{Construction}

\newtheorem{remark}[definition]{Remark}
\newenvironment{rema}{\begin{remark} \rm}{\end{remark}}
\newtheorem{remarks}[definition]{Remarks}

\newtheorem{example}[definition]{Example}
\newenvironment{exam}{\begin{example} \rm}{\end{example}}
\newtheorem{examples}[definition]{Examples}

\newtheorem{nothing}[definition]{$\!\!$}

\newenvironment{proo}{{\flushleft \it Proof.}}{\hfill $\square$
  \vspace{2mm}}
\newenvironment{proo-theo}{{\flushleft \it Proof of Theorem
    \ref{theo-transplanting}.}}{\hfill $\square$ \vspace{2mm}}

\newtheorem{conjecture}[definition]{Conjecture}

\newtheorem{definition*}{Definition}[section]
\newenvironment{defi*}{\begin{definition*} \rm}{\end{definition*}}
\newtheorem{definitions*}[definition*]{Definitions}
\newenvironment{defis*}{\begin{definitions*} \rm}{\end{definitions*}}
\newtheorem{prop*}[definition*]{Proposition}
\newtheorem{lemm*}[definition*]{Lemma}
\newtheorem{coro*}[definition*]{Corollary}
\newtheorem{theo*}[definition*]{Theorem}
\newtheorem{remark*}[definition*]{Remark}
\newenvironment{rema*}{\begin{remark*} \rm}{\end{remark*}}
\newtheorem{remarks*}[definition*]{Remarks}
\newenvironment{remas*}{\begin{remarks*} \rm}{\end{remarks*}}
\newtheorem{example*}[definition*]{Example}
\newenvironment{exam*}{\begin{example*} \rm}{\end{example*}}
\newtheorem{examples*}[definition*]{Examples}
\newenvironment{exams*}{\begin{examples*} \rm}{\end{examples*}}
\newtheorem{nothing*}[definition*]{$\!\!$}
\newenvironment{noth*}{\begin{nothing*} \rm}{\end{nothing*}}

\newtheorem{commentaire*}[definition*]{Commentaire}

\def \pic {{{\rm Pic}}}
\def \codim {{\rm Codim}}
\def \sca #1#2{\langle #1,#2 \rangle}
\def\s{\sigma}
\newcommand{\G}{\mathbb{G}}

\newcommand{\Z}{\mathbb Z}

\newcommand{\Q}{\mathbb{Q}}
\newcommand{\C}{\mathbb{C}}
\renewcommand{\O}{\mathbb{O}}
\newcommand{\F}{\mathbb{F}}
\newcommand{\p}{\mathbb{P}}

\renewcommand{\a}{{\alpha}}

\newcommand{\pu}{{\mathbb{P}^1}}

\begin{document}

\title{Small codimension subvarieties in homogeneous spaces}
\author{N. Perrin}
\date{}

\maketitle

\begin{abstract}
We prove Bertini type theorems for the inverse image, under a proper
morphism, of any Schubert variety in an homogeneous space. Using generalisations of Deligne's trick, we deduce connectedness results for the inverse image of the diagonal in $X^2$ where $X$ is any isotropic grassmannian. We also deduce simple connectedness properties for subvarieties of $X$. Finally we prove transplanting theorems {\`a} la Barth-Larsen for the Picard group of any isotropic grassmannian of lines and for the Neron-Severi group of some adjoint and coadjoint homogeneous spaces.
\end{abstract}

 {\def\thefootnote{\relax}
 \footnote{ \hspace{-6.8mm}
 }
 }

\begin{center}{\bf Introduction}\end{center}

In this text we work over an algebraically closed field of 
characteristic zero.
The topological properties of subvarieties of small codimension in a
rational homogeneous space $X$ have drawn much attention since the
work of W. Barth and M.E. Larsen \cite{barth-larsen} on the
projective space. The litterature is very large on the subject, let us
only mention the celebrated paper of R. Hartshorne \cite{hartshorne}
where the so called Harshorne's conjecture for small codimension
subvarieties of the projective space was formulated and where he gives
a simple proof of Barth-Larsen's theorem using the hard Lefschetz
Theorem. We refer to Lazarsfeld's book \cite[Chapter 3]{lazarsfeld}
for a more extensive review of the litterature on these themes. 

Some of
the most achieved results on this subject can be stated using the
following invariant. Recall that for $X$ a rational homogeneous space, the tangent bundle $T_X$ of $X$ is globally generated.

\begin{defi}
  \label{defi-coampleness}
 Let $\phi:\p(T_X^\vee)\to\p^N$ be the map defined by the global sections of the tangent bundle $T_X$. The \emph{ampleness} of $X$, denoted by ${\rm a}(X)$ is the maximum of the dimensions of the fibres of $\phi$. The \emph{coampleness} of $X$ denoted by ${\rm ca}(X)$ is $\dim X-{\rm a}(X)$.
\end{defi}

G. Faltings in \cite{faltings} proves the following theorem.

\begin{theo}[Faltings]
  \label{theo-faltings}
Let $X$ be a rational homogeneous space and $f:Z\to X\times X$ a proper map with $Z$ irreducible and $\codim Z< {\rm ca}(X)$, then $f^{-1}(\Delta)$ is connected where $\Delta$ is the diagonal in $X\times X$.
\end{theo}

As a corollary, one obtains that any irreducible subvariety $Y$ in $X$ satisfying the inequality $2\codim Y\leq {\rm ca}(X)-1$ is algebraically simply connected \emph{i.e.} $\pi_1^{\rm alg}(Y)=1$. The invariant ${\rm ca}(X)$ also appears in the work of A. Sommese. For example, the following theorem is a consequence of results he proves with A. van de Ven in \cite{sommese-vandeven}.

\begin{theo}[Sommese-van de Ven]
  \label{theo-sommese-van-de-ven}
Let $X$ be a rational homogeneous space and $Y$ be a smooth connected subvariety of $X$, then for any point $y$ in $Y$ we have the vanishing $\pi_j(X,Y,y)=0$ for $2\codim Y\leq {\rm ca}(X)-j+1$.
\end{theo}

In particular, they recover the above (algebraically) simple connectedness for smooth subvarieties and furthermore, if $2\codim Y\leq {\rm ca}(X)-2$, then they obtain the equality $\pic(Y)=\pic(X)$. These results are made very useful by the computation of the invariant ${\rm ca}(X)$ for any rational homogeneous space by N. Goldstein \cite{goldstein}. 

\vskip 0.5 cm

In this paper, we shall only deal with homogeneous spaces with Picard number 1. We want to study the properties, in rational homogeneous spaces, of small codimension subvarieties having additional numerical properties but higher codimension. The fact that numerical properties of $Y$ should be sufficient in this setting was suggested by W. Fulton and R. Lazarsfeld in \cite{fulton-lazarsfeld} and verified for product of projective spaces and grassmannians by O. Debarre in \cite{debarre}. In particular, if $X$ is the grassmannian $\G(p,n)$ of $p$-dimensional subspaces in an $n$-dimensional vector space and if $\s$ (resp. $\s'$) is the Schubert class of subspaces containing a fixed vector (resp. contained in a fixed hyperplane), then O. Debarre proves the following result (for $Y\subset X$, we denote by $[Y]$ the cohomology class of $Y$).

\begin{theo}[Debarre]
  \label{theo-debarre}
Let $f:Z\to X\times X$ be a morphism with $Z$ irreducible and denote by $i:X\to X\times X$ the diagonal embedding and by $\Delta$ the diagonal. If $[f(Z)](i_*\s+i_*\s')\neq0$, then for $\Delta^g$ a general translate of $\Delta$, we have that $f^{-1}(\Delta^g)$ is connected.

In particular, if $Y$ is an irreducible subvariety of $X$ with $[Y]^2(\s+\s')\neq0$, then $Y$ is simply connected.
\end{theo}

This result is obtained with the same method as the one used in \cite{fulton-lazarsfeld} for the projective space. First Debarre proves a Bertini type Theorem for inverse images of Schubert varieties in the grassmannian. Then by adapting to the situation a trick of Deligne, he is able to replace the inverse image of the diagonal in $X\times X$ by the inverse image of a Schubert variety in a larger grassmannian.

\vskip 0.5 cm

In the first section of this paper, we prove Bertini type Theorems for inverse images of Schubert varieties under proper morphisms. We define for any Schubert variety $X^P(w)$ in $X$ the notion of an \emph{admissible} Schubert subvariety (see Definition \ref{defi-admissible}) and prove the following theorem.

\begin{theo}
\label{theo-main1}
  Let $X$ be an homogeneous space, let $f:Y\to X$ be a proper map with $Y$ irreducible. Assume that for some admissible Schubert subvariety $X^P(v)$ of a Schubert variety $X^P(w)$ we have $[f(Y)][X^P(v)]\neq0$, then $f^{-1}(X^P(w))$ is connected.
\end{theo}

In the second section, we apply the above theorem to
obtain results on the connectedness of the inverse image of the diagonal. For this we need
to further generalise Deligne's trick. We succeded to generalise it
only for classical groups therefore we only obtain in this section
results for homogeneous spaces under classical groups. More precisely,
let $\G_Q(p,n)$ and $\G_\omega(p,2n)$ be the grassmannians of
isotropic subspaces of dimension $p$ in a vector space (of dimension $n$ resp. $2n$) endowed with a
non degenerate quadratic form $Q$ resp. symplectic form $\omega$. Let
$\s$ be the Schubert class in $\G_Q(p,n)$ resp. $\G_\omega(p,2n)$ of
subspaces contained in the orthogonal of an isotropic vector, we
obtain the following theorem. 

\begin{theo}
\label{theo-main2}
  Let $f:Z\to X\times X$ be a morphism with $Z$ irreducible and denote by $i:X\to X\times X$ the diagonal embedding and by $\Delta$ the diagonal. If $[f(Z)]i_*\s\neq0$, then $f^{-1}(\Delta)$ is connected.

In particular, if $Y$ is an irreducible subvariety of $X$ with $[Y]^2\s\neq0$, then $Y$ is simply connected.
\end{theo}

Note that using the first section, we can prove the above result only for proper morphisms and therefore we only obtain algebraically simple connectedness of $Y$. In the last section we explain how to remove the proper hypothesis.

In the third section, we prove, using the first two sections, transplanting theorems for the group of divisors modulo numerical equivalence. This was our first motivation for studying connectedness properties. Such a transplanting theorem was proved by E. Arrondo and J. Caravantes in \cite{arrondo-caravantes} for the Picard group of subspaces in the grassmannian $\G(2,n)$ only using numerical properties and the connectedness results of Debarre \cite{debarre}.

\begin{defi}
\label{defi-cumbersome}
  A subvariety $Y$ of a variety $X$ is called \emph{cumbersome} if for all subvariety $Z$ of $X$ with $\dim Z=\codim Y$, we have $[Y][Z]\neq0$.
\end{defi}

\begin{theo}[Arrondo-Caravantes]
  Let $Y$ be a smooth cumbersome subvariety of $\G(2,n)$ with $2\codim Y\leq\dim X-2$, then $\pic(Y)=\Z$.
\end{theo}

The proof of this theorem involves a mysterious (as the author themselves confess) computation leading to expressing some quadratic form as sum of squares. In \cite{cras}, we explained that their method is related to the fact that the intersection form on the middle cohomology group is positive definite. We also explained how the method combined with connectedness result would also apply to any other homogeneous space with the same property. In the third section we revisit this method and make it work in a more general context. We shall need the definitions.

\begin{defi}
\label{defi-eff}
  (\i) Let $X$ be a rational homogeneous space with Picard number 1 and let $h$ be an ample generator of the Picard group. We define the bilinear form $(\ ,\ )_{h^k}$ on $H^{\dim X-k}(X,\Z)$ by $(\s,\s')_{h^k}=\s\s'h^k$.

(\i\i) We define the \emph{effectiveness} of $X$, by ${\rm eff}(X)=\min\{k\ /\ (\ ,\ )_{h^k} \textrm{ is positive definite}\}$. The \emph{coeffectiveness} of $X$ is ${\rm coeff}(X)=\dim X-{\rm eff}(X)$.
\end{defi}

Denote by $N^1(X)$ is the group of divisors modulo numerical equivalence in $X$. The main result of the third section is

\begin{theo}
\label{theo-main3}
  Let $X$ be a rational homogeneous space with Picard number 1. Let $Y$ be a smooth cumbersome subvariety in $X$. If $2\codim Y\leq {\rm coeff}(X)-2$, then $N^1(X)=\Z$.
\end{theo}

This theorem should be though as a numerical version of Theorem \ref{theo-sommese-van-de-ven} for divisors. However, this theorem does not improve Sommese and van de Ven theorem for all homogeneous spaces. This is the case when ${\rm coeff}(X)>{\rm ca}(X)$ and only occurs for projective spaces over composition algebras: $\p^n$, $\G(2,n)$ and $\mathbb{O}\p^2$, for adjoint varieties (also called minimal adjoint orbits) and coadjoint varieties, see Definition \ref{defi-adjoint}. As a consequence of this result and of the simple connectedness properties obtained in the second section we have

\begin{coro}
  Let $X$ be $\G_Q(2,2n+1)$, $\G_\omega(2,2n)$ or $\G_Q(2,4n)$ and $Y$ be a smooth cumbersome subvariety with $2\codim Y\leq{\dim X-3}$, then $\pic(Y)=\Z$.
\end{coro}

In the last section, we extend the results of the second section for any (non necessarily proper) map and deduce results for the topological fundamental group instead of results for the algebraic fundamental group. This extension is not needed for the transplanting results of the third section.

\tableofcontents

\section*{Notation and conventions}

We will work over an algebraically closed field of characteristic zero. We shall follow the notation from \cite{fulton-lazarsfeld} and \cite{debarre}. In particular, we shall forget the base point in the notation $\pi_1(Y)$ of the fundamental group of $Y$. We shall also at several occasion not repeat the proofs of \cite{fulton-lazarsfeld} and refer to this text for the arguments instead of recopying them. 

\vskip 0.2 cm

We fix $G$ a semisimple algebraic group, let us also fix $T\subset B\subset P\subset G$ a maximal torus, a Borel subgroup and a maximal parabolic subgroup in $G$. We denote by $X$ the rational homogeneous space $G/P$. Note that we have $\pic(X)=\Z$. 
We denote by $W$ the Weyl group of $G$ and by $W_P$ the Weyl group of $P$. The Schubert varieties in $G/P$ (\emph{i.e.} the closures of the $B$-orbits) are indexed by the set $W^P$ of minimal length representatives of the quotient $W/W_P$. For any element $w$ in $W^P$, we denote by $X^P(w)$ the corresponding Schubert variety and by $\s(w)\in H^{2l(w)}(X,\Z)$ the corresponding cohomology class with $l(w)$ the length of $w$. For $\a$ a root of $G$, we denote by $U(\a)$ the corresponding unipotent subgroup in $G$. For $Q$ a parabolic subgroup containing $B$, we denote by $\Sigma(Q)$ the set of simple roots such that $U(-\a)$ is not contained in $Q$. We refer to the text of M. Brion \cite{brion} for more details on the geometry of rational homogeneous spaces.

We shall also use the following less usual notation. We denote by $S^P(w)$ the subgroup of $G$ stabilising the Schubert variety $X^P(w)$. The subgroup $S^P(w)$ is a parabolic subgroup of $G$ containing $B$ and we set $\Sigma^P(w)=\Sigma(S^P(w))$. 

\vskip 0.2 cm

\paragraph{(Co)adjoint homgeneous spaces.} We shall now define the adjoint and coadjoint varieties of the group $G$. Let us denote by $\Theta$ (resp. $\theta$) the highest root of $G$ (resp. the highest short root of $G$). If $G$ is simply laced then $\Theta=\theta$. For $\varpi$ a dominant weight, we denote by $V_\varpi$ highest weight representation of $G$ with highest weight $\varpi$ and by $v_\varpi$ an highest weight vector.

\begin{defi}
  \label{defi-adjoint}
The \emph{adjoint} (resp. \emph{coadjoint}) variety of $G$ is the (closed) orbit in $\p(V_\varpi)$ of $[v_\varpi]$ under the action of $G$ with $\varpi=\Theta$ (resp. $\varpi=\theta$).
\end{defi}

\section{Bertini type theorems}
\label{section-bertini}

In this section we shall prove Theorem \ref{theo-main1} on the irreducibility of inverse images of Schubert varieties under proper maps from irreducible sources.

\subsection{Bertini for minimal generating Schubert varieties}

\begin{defi}
  Let us call \emph{minimal generating Schubert variety} of $X$ any Schubert subvariety $X^P(t)$ for $t\in W^P$ such that $t$ has a reduced expression $s_{\beta_1}\cdots s_{\beta_k}$ with all the simple roots $\beta_i$ distinct and $\{\beta_1,\cdots,\beta_k\}=\Sigma(P)$.
\end{defi}

\begin{rema}
 (\i) We use the terminology of generating Schubert varieties because
 these Schubert varieties are generating subvarieties in the sense of
 Chow \cite{chow} and therefore are G3 in $X$ (see also \cite{badescu}
 for more details on the formal geometry point of view). 

(\i\i) The minimal generating Schubert varieties are the smallest Schubert varieties with positive degree with respect to any non trivial element of the monoid of pseudo-effective divisors.
\end{rema}

\begin{prop}
\label{prop-bertini-min}
  Let $Y$ be an irreducible variety, let $f:Y\to X$ be a proper dominant morphism and let $X^P(t)$ be a minimal generating Schubert variety in $X$.

(\i) The inverse image $f^{-1}(g\cdot X^P(t))$ is irreducible for $g$ in a dense open subset of $G$.

(\i\i) The inverse image $f^{-1}(g\cdot X^P(t))$ is connected for all $g$ in $G$.

(\i\i\i) If $Y$ is unibranch, then $\pi_1^{\rm alg}(f^{-1}(g\cdot
X^P(t)))\to\pi_1^{\rm alg}(Y)$ is surjective for $g$ general in $G$.
\end{prop}

\begin{proo}
  We use Proposition 1 in \cite{paranjape-srinivas} proving \textit{(\i)} and \textit{(\i\i)} for $f$ a finite map. By Stein factorisation, we obtain \textit{(\i\i)} for any proper map and by Kleiman-Bertini Theorem \cite{kleiman} we know that for $g\in G$ general the inverse image $f^{-1}(g\cdot X^P(t))$ is locally integral (for Zariski topology) and because it is connected it has to be irreducible.

We get \textit{(\i\i\i)} by applying \textit{(\i)} to $\widetilde{f}=f\circ\pi$ where $\pi:\widetilde{Y}\to Y$ is any connected {\'e}tale covering. Note that the {\'e}tale covering being connected and locally irreducible (because $Y$ is unibranch), it is irreducible. We get a connected {\'e}tale covering $\widetilde{f}^{-1}(g\cdot X^P(w))\to f^{-1}(g\cdot X^P(w))$ by restriction of the previous one. This proves the result.
\end{proo}

\subsection{Bertini for Schubert varieties}

\begin{defi}
\label{defi-admissible}
  Let $X^P(v)\subset X^P(w)$ be an inclusion of Schubert varieties. We shall say that $X^P(v)$ is \emph{admissible} in $X^P(w)$ if we have $S^P(w)\cdot X^P(v)=X^P(w)$ and $\Sigma^P(w)\cap \Sigma^P(v)=\emptyset$.
\end{defi}

\begin{theo}
\label{theo-bertini}
  Let $f:Y\to X$ be a proper morphism with $Y$ irreducible and such that there exists an admissible Schubert subvariety $X^P(v)$ of the Schubert variety $X^P(w)$ with $[f(Y)]\cdot[X^P(v)]\neq0$. 

(\i) Then the inverse image $f^{-1}(g\cdot X^P(w))$ is irreducible for $g$ in an open subset of $G$.

(\i\i) Then the inverse image $f^{-1}(g\cdot X^P(w))$ is connected for all $g$ in $G$.

(\i\i\i) If $Y$ is unibranch, then $\pi_1^{\rm alg}(f^{-1}(g\cdot X^P(w)))\to\pi_1^{\rm alg}(Y)$ is surjective for $g$ general in $G$.

(\i v) If $Y$ is unibranch, for any $g\in G$ and any non trivial neighbourhood $U$
of $g\cdot X^P(w)$, the map $\pi_1^{\rm alg}(f^{-1}(U))\to\pi_1^{\rm
  alg}(Y)$ is surjective.
\end{theo}

\begin{proo}
We deduce \textit{(\i\i)} and \textit{(\i\i\i)} from \textit{(\i)} as in Proposition \ref{prop-bertini-min}. Let us prove \textit{(\i)}.

  Let $Q=S^P(w)$ be the stabiliser of $X^P(w)$ and let us consider the following subvariety 
$I^P_Q(w)=\{(\bar g,\bar h)\in G/P\times G/Q\ /\ \bar g\in h\cdot X^P(w)\}$
of the product $G/P\times G/Q$.
This is well defined since $Q$ stabilises $X^P(w)$ and we remark that
the two projection $p$ resp. $q$ to $G/P$ resp. $G/Q$ realise
$I_Q^P(w)$ as a locally trivial fibration with fiber isomorphic to
$X^Q(u)$ resp. $X^P(w)$ where the Schubert variety $X^Q(u)$ is the
closure in $G/Q$ of the $P$-orbit of the Schubert variety
$X^Q(w^{-1})$. In particular, the fiber product
$Z=Y\times_{G/P}I_Q^P(w)$ is irreducible as a locally trivial
fibration with fiber $X^Q(u)$ over $Y$. Because of our assumption, the image
of $f$ meets general translates of $X^P(w)$ therefore the composition
$f':Z\to I_Q^P(w)\to G/Q$ is proper and dominant. 

  \begin{lemm}
\label{lemm-section}
    For any minimal generating Schubert variety $X^Q(t)$ in $G/Q$ and
    for any $y\in Y$ with $f(y)\in X^P(v)$, the map $\bar h\mapsto
    (y,\bar h)$ defines a section $s:X^Q(t)\to Z$ of $f'$.
  \end{lemm}

  \begin{proo}
    We only need to prove that for any element $\bar h\in X^Q(t)$, the
    element $(f(y),\bar h)$ lies in $I^P_Q(w)$. Let us first remark
    that if $t=s_{\beta_1}\cdots s_{\beta_k}$ is a reduced expression
    of $t$, then the Schubert variety $X^Q(t)$ is the closure of
    $U(\beta_1)\cdots U(\beta_k)\cdot \bar e$ where $U(\beta)$ is the
    unipotent subgroup associated to the root $\beta$ and $e$ is the
    unit element in $G$. Let us also remark that for $f(y)\in
    X^P(v)$, because of the inclusion $X^P(v)\subset X^P(w)$, we have
    $(f(y),\bar e)\in I_Q^P(w)$. Now let $U(\beta_k)\cdots
    U(\beta_1)$ act on the inclusions $f(y)\in X^P(v)\subset
    X^P(w)$. We get, by definition of a minimal generating Schubert
    variety and because $X^P(v)$ is admissible the inclusions
    $U(\beta_k)\cdots U(\beta_1)\cdot f(y)\subset X^P(v)\subset
    X^P(w)$. This is equivalent to the inclusion $f(y)\in
    U(\beta_1)\cdots U(\beta_k)\cdot X^P(w)$ and gives the inclusion
    $(f(y),U(\beta_1)\cdots U(\beta_k)\cdot\bar e)\in I_Q^P(w)$,
    proving the result.
  \end{proo}

We now want to apply the following result (see for example \cite[Lemma
3.3.2]{lazarsfeld}). 

\begin{lemm}
\label{lemm-lazarsfeld}
  Let $q:V\to T$ be a dominating morphism of irreducible complex
  varieties such that $q$ admits a section $s:T\to V$ whose image does
  not lie in the singular locus of $V$. Then a general fiber of $q$ is
  irreducible.
\end{lemm}

Let us prove some regularity results. First remark that because the map $I_Q^P(w)$ is a locally trivial fibration with fiber $X^P(w)$ over $G/Q$, then the smooth locus $I^P_Q(w)^{\rm sm}$ of $I_Q^P(w)$ is the associated locally trivial fibration over $G/Q$ with fiber $X^P(w)^{\rm sm}$ the smooth locus of $X^P(w)$. Let $Y^{\rm sm}$ be the smooth locus of $Y$, then the fiber product $Y^{\rm sm}\times_{G/P}I_Q^P(w)^{\rm sm}$ is smooth. If $\Omega^Q(t)$ is the Schubert cell in $X^Q(t)$, by Kleiman-Bertini theorem (see \cite{kleiman}), there exists a dense open subset $U$ of $G$ such that for $g\in U$, the fiber product $(Y^{\rm sm}\times_{G/P}I_Q^P(w)^{\rm sm})\times_{G/Q}g\cdot \Omega^Q(t)$ is smooth. Note that, we have the equality 
$$(Y^{\rm sm}\times_{G/P}I_Q^P(w)^{\rm sm})\times_{G/Q}g\cdot \Omega^Q(t)=\{(y,\bar h)\ /\ \bar h\in\Omega^Q(t),\ y\in Y^{\rm sm}\ {\rm and}\ f(y)\in gh\cdot X^P(w)^{\rm sm}\}.$$
As we have $[f(Y)]\cdot[X^P(v)]\neq0$, restricting $U$ and using Kleiman-Bertini again, we may assume that there exists $y\in Y^{\rm sm}$ with $f(y)\in g\cdot \Omega^P(v)$. By assumption $X^P(v)$ is admissible thus $\Omega^P(v)$ lies in $X^P(w)^{\rm sm}$, we get that $(y,\bar e)$ lies in $(Y^{\rm sm}\times_{G/P}I_Q^P(w)^{\rm sm})\times_{G/Q}g\cdot \Omega^Q(t)$. Therefore, by Lemma \ref{lemm-section}, we have the inclusion $(y,\bar h)\in (Y^{\rm sm}\times_{G/P}I_Q^P(w)^{\rm sm})\times_{G/Q}g\cdot \Omega^Q(t)$ for $\bar h$ in a dense open subset of $\Omega^Q(t)$. Now remark that $Y^{\rm sm}\times_{G/P}I_Q^P(w)^{\rm sm}$ is a subvariety of $Z=Y\times_{G/P}I^P_Q(w)$ and that $(Y^{\rm sm}\times_{G/P}I_Q^P(w)^{\rm sm})\times_{G/Q}g\cdot \Omega^Q(t)$ is a subvariety of $Z\times_{G/Q}X^Q(t)$. The map $f':Z\to G/Q$ is dominating thus by Proposition \ref{prop-bertini-min}, we know that, restricting $U$ if necessary, the variety $Z\times_{G/Q}g\cdot X^Q(t)$ is irreducible. The same is true for $Z\times_{G/Q}g\cdot\Omega^Q(t)$. By the above we have a section in the smooth locus, thus by Lemma \ref{lemm-lazarsfeld}, the general fiber of $f'$ over $\bar h\in g\cdot\Omega^Q(t)$ is irreducible. But this fiber is exactely 
$${f'}^{-1}(\bar h)=\{y\in Y\ /\ f(y)\in gh\cdot
X^P(w)\}=f^{-1}(gh\cdot X^P(w))$$
and the result follows.

For \textit{(\i v)}, any $U$ contains
translates $g'\cdot X^P(w)$ for which \textit{(\i\i\i)} apply.
\end{proo}

As an easy corollary of this result, we partially answer a question raised in \cite{cras}. Let $X$ be $\O\p^2(\C)$ the Cayley plane. It is the homogeneous space $G/P$ with $G$ a group of type $E_6$ and $P$ a maximal parabolic subgroup with $\Sigma(P)=\{\a_1\}$ with $\a_1$ a simple root with notation as in \cite{bou}.

\begin{coro}
  Let $Y$ be a smooth subvariety of $X$. 

(\i) If $\dim Y\geq12$, then $\pic(Y)=\Z$.

(\i\i) If $\dim Y\geq11$ and $Y$ is cumbersome, then $\pic(Y)=\Z$.

(\i\i\i) If $\dim Y\geq9$ and $Y$ is cumbersome, then $N^1(Y)=\Z$. If furthermore $Y$ is simply connected, then $\pic(Y)=\Z$.
\end{coro}

\begin{proo}
  For \textit{(\i)}, we only apply the result of Sommese and van de Ven, see Theorem \ref{theo-sommese-van-de-ven}, and the fact proved by Goldstein \cite{goldstein} that in this case ${\rm ca}(X)=11$.

For \textit{(\i\i)} and \textit{(\i\i\i)}, the results in \cite{cras} imply that if $Y$ is cumbersome of dimension at least 9 such that the intersection of $Y$ with a general translate of a dimension 9 Schubert variety is irreducible, then $N^1(Y)=\Z$. If furthermore $Y$ is simply connected, we get that $\pic(Y)=\Z$ (to get the simply connectedness in \textit{(\i\i)} we again apply the results of Sommese and Goldstein). Therefore we only need to prove the irreducibility of the intersection of $Y$ with a general translate of a dimension 9 Schubert variety. There are only two such Schubert varieties $X^P(w)$ whose reduced expressions are $w=s_1s_3s_4s_2s_6s_5s_4s_3s_1$ resp. $w=s_5s_3s_4s_2s_6s_5s_4s_3s_1$ (here the notation are those of \cite{bou}). Then the Schubert variety $X^P(v)$ with $v$ given by $v=s_1s_3s_4s_2s_5s_4s_3s_1$ resp. $v=s_5s_4s_2s_6s_5s_4s_3s_1$ are admissible in $X^P(w)$ and because $Y$ is cumbersome of dimension at leat 9, we have $[f(Y)]\cdot[X^P(v)]\neq0$ and the result follows from the previous theorem.
\end{proo}

\section{Connectedness Theorems}
\label{section-connexite}

To prove connectedness theorems for the inverse image of the diagonal, we will use our results on Schubert varieties together with Deligne's trick as explained in \cite{fulton-lazarsfeld} and \cite{lazarsfeld} for the projective space or in \cite{debarre} for the grassmannian variety. However, for exceptional groups, we were not able to find a suitable generalisation of Deligne's trick, so our results here are valid only for classical groups.

\subsection{Deligne's trick}

In this section, we generalise, following and generalising O. Debarre \cite{debarre}, Deligne's trick to reduce the connectedness of the inverse image of the diagonal to Bertini type theorems. 

Let $V$ be a vector space endowed with a non degenerate quadratic form $Q$ (resp. symplectic form $\omega$). Let us consider $W=V_1\oplus V_2$ where both $V_1$ and $V_2$ are isomorphic to $V$. Define over $W$ a non degenerate quadratic (resp. symplectic) form by looking at $W$ as the orthogonal sum of $V_1$ and $V_2$. For $i\in\{1,2\}$ and for $v_i\in V_i$, we denote by $[v_i]$ and $[v_1,v_2]$ the classes of $v_i$ and $(v_1,v_2)$ in $\p(V_i)$ and $\p(W)$. 

We denote by $\G_Q(p,W)$ (resp. $\G_\omega(p,W)$) the grassmannian of isotropic $p$-dimensional vector subspaces of $W$. We shall denote by $\G_Q^0(p,W)$ (resp. $\G_\omega^0(p,W)$), the open subvarieties defined by $\{V_p\ /\ V_p\cap V_i=0 \textrm{ for } i\in\{1,2\}\}$.

\begin{lemm}
\label{lemm-deligne}
 Let $p_Q$ (resp. $p_\omega$) denote the map $\G_Q^0(p,W)\to \G(p,V)\times \G(p,V)$ (resp. the map $\G_\omega^0(p,W)\to \G(p,V)\times \G(p,V)$) defined by the two projections from $W$ to $V_i$ for $i\in\{1,2\}$. Then the restriction of $p_Q$ (resp. $p_\omega$) over $\G_Q(p,V)\times \G_Q(p,V)$ (resp. over $\G_\omega(p,V)\times \G_\omega(p,V)$) is a ${\rm GL}_p(\C)$-bundle.
\end{lemm}

\begin{proo}
 The choice of a subspace  in the fiber $W_p$ over a pair $(V_p,V_p')$ of subspaces determines a linear isomorphism $\gamma:V_p\to V_p'$: the graph of $\gamma$ is the subspace $W_p$. Now if the two subspaces $V_p$ and $V_p'$ are isotropic, then there is no condition on $\gamma$ for $W_p$ to be isotropic and the result follows.
\end{proo}

\subsection{Inverse image of the diagonal}

Let $v$ be an isotropic vector in $V$ and define a Schubert variety in $\G_Q(p,V)$ (resp. $\G_\omega(p,V)$) by $\Gamma=\{V_p\ /\ V_p\subset v^\perp\}$.
We denote by $\Delta$ the diagonal in the product $\G_Q(p,V)\times\G_Q(p,V)$ (resp. $\G_\omega(p,V)\times\G_\omega(p,V)$). We embed $\Gamma$ in $\Delta$ and then in the product $\G_Q(p,V)\times\G_Q(p,V)$ (resp. $\G_\omega(p,V)\times\G_\omega(p,V)$). Finally, we denote by $\gamma$ the cohomology class of this embedded variety. This class is of degree $\dim\G_Q(p,V)+\dim V-p-1$ (resp. $\dim\G_\omega(p,V)+\dim V-p$).

\begin{theo}
\label{theo-connexite}
  Let $f:Z\to \G_Q(p,V)\times\G_Q(p,V)$ (resp. $f:Z\to \G_\omega(p,V)\times\G_\omega(p,V)$) be a proper morphism with $Z$ irreducible. Assume that $[f(Z)]\cdot\gamma\neq0$.

(\i) Then $f^{-1}(\Delta)$ is connected.

(\i\i) If $Z$ is unibranch, then $\pi_1^{\rm alg}(f^{-1}(\Delta))\to\pi_1^{\rm alg}(Z)$ is surjective.
\end{theo}

\begin{proo}
  \textit{(\i)} Let us consider the ${\rm GL}_p(\C)$-bundle $p_Q$ (resp. $p_\omega$) over the product $\G_Q(p,V)\times\G_Q(p,V)$ (resp. $\G_\omega(p,V)\times\G_\omega(p,V)$) as decribed in Lemma \ref{lemm-deligne}. The fiber product $Z'$ constructed from this bundle and $f$ is again irreducible and mapped to $\G_Q(p,W)$ (resp. $\G_\omega(p,W)$). 

Let us denote by $\Delta_W$ the diagonal embedding of $V$ in $W=V\oplus V$. Let $W_1$ and $W_2$ be general isotropic subspaces of dimension $\lfloor\dim V/2\rfloor$ in $\Delta_W^\perp$. If $\dim V$ is odd, we also fix a general isotropic line $L$ in $\Delta_W^\perp$ (if $\dim V$ is even we set $L=0$ to simplify notation). Note that $W_1$, $W_2$ and $L$ are in direct sum and that $W_1\oplus W_2\oplus L=\Delta_W^\perp$ thus $(W_1\oplus W_2\oplus L)^\perp=\Delta_W$. Therefore, if we define $\Delta\G=\{V_p\ /\ V_p\subset \Delta_W\}$, we have the following equality
$$\Delta\G=\{V_p\ /\ V_p\subset W_1^\perp\}\cap \{V_p\ /\ V_p\subset W_2^\perp\}\cap \{V_p\ /\ V_p\subset L^\perp\}.$$
Moreover, the projection $p_Q$ (resp. $p_\omega$) realises an isomorphism from $\Delta\G$ onto $\Delta$ thus we have an isomorphism between ${f'}^{-1}(\Delta\G)$ and $f^{-1}(\Delta)$. We therefore need to prove the connectedness of ${f'}^{-1}(\Delta\G)$. But the three varieties $\{V_p\ /\ V_p\subset W_1^\perp\}$, $\{V_p\ /\ V_p\subset W_2^\perp\}$ and $\{V_p\ /\ V_p\subset L^\perp\}$ are Schubert varieties and we can apply the results from Section \ref{section-bertini}. First remark that if $W'_i$ and $L'$ are isotropic subspaces in $W$ such that $W_i$ is an hyperplane in $W'_i$ (resp. $L$ is an hyperplane in $L'$), then the Schubert varieties $\{V_p\ /\ V_p\subset {W_1}'^\perp\}$, $\{V_p\ /\ V_p\subset {W'_2}^\perp\}$ and $\{V_p\ /\ V_p\subset {L'}^\perp\}$ are admissible Schubert subvarieties of $\{V_p\ /\ V_p\subset W_1^\perp\}$, $\{V_p\ /\ V_p\subset W_2^\perp\}$ and $\{V_p\ /\ V_p\subset L^\perp\}$. 
By applying Theorem \ref{theo-bertini} three times, the variety ${f'}^{-1}(\Delta\G)$ is connected as soon as the product of $[f'(Z')]$ with $[\{V_p\ /\ V_p\subset {W_1'}^\perp\}]$, with $[\{V_p\ /\ V_p\subset {W_1}^\perp\}]\cdot[\{V_p\ /\ V_p\subset {W_2}^\perp\}]$ and with $[\{V_p\ /\ V_p\subset {W_1}^\perp\}]\cdot[\{V_p\ /\ V_p\subset {W_2}^\perp\}]\cdot[\{V_p\ /\ V_p\subset {L'}^\perp\}]$ are non zero. But these classes are described as $[\{V_p\ /\ V_p\subset {W_1'}^\perp\}]$, $[\{V_p\ /\ V_p\subset (W_1\oplus W'_2)^\perp\}]$ and $[\{V_p\ /\ V_p\subset (W_1\oplus W_2\oplus L')^\perp\}]$, thus it is clear that the non vanishing of the product with the last one implies the two others. Therefore it is enough to prove that the multiplication of $[f(Z)]$ with the class of the image via $p_Q$ (resp. $p_\omega$) of $\{V_p\ /\ V_p\subset (W_1\oplus W_2\oplus L')^\perp\}$ is non trivial. As the class of this last variety is equal to $\gamma$, the result follows.

For \textit{(\i\i)}, we apply the standard tricks as explained in \cite{fulton-lazarsfeld} in Remark 2.2 and page 40.
\end{proo}

Let $v$ be an isotropic point in $V$ and let us denote by $\s$ the cohomology class of the Schubert variety $\{V_p\ /\ V_p\subset v^\perp\}$. 

\begin{coro}
  Let $f:Y\to\G_Q(p,V)$ and $g:Y'\to\G_Q(p,V)$ (resp. $f:Y\to\G_\omega(p,V)$ and $g:Y'\to\G_\omega(p,V)$) be proper morphisms such that $[f(Y)]\cdot [g(Y')]\cdot \s\neq0$.

(\i) Then $Y\times_{\G_Q(p,V)}Y'$ (resp. $Y\times_{\G_\omega(p,V)}Y'$) is connected.

(\i\i) If the varieties $Y$ and $Y'$ are unibranch, then the map $\pi_1^{\rm alg}(Y\times_{\G_Q(p,V)}Y')\to\pi_1^{\rm alg}(Y\times Y')$ (resp. $\pi_1^{\rm alg}(Y\times_{\G_\omega(p,V)}Y')\to\pi_1^{\rm alg}(Y\times Y')$) is surjective.
\end{coro}

\begin{rema}
  We will prove in the last section that the above theorem and corollary hold true witout the assumption $f$ proper. We therefore obtain stronger results with $\pi_1$ in place of $\pi_1^{\rm alg}$.
\end{rema}

\begin{proo}
  We apply the former theorem to $f\times g$. We only need to prove that $[f(Y)\times g(Y')]\cdot\gamma\neq0$. But if $i$ is the closed immersion of $\Delta$ in $\G_Q(p,V)\times\G_Q(p,V)$ (resp. $\G_\omega(p,V)\times\G_\omega(p,V)$), we have $\gamma=i_*\s$ thus we have the equalities $[f(Y)\times g(Y')]\cdot\gamma=i_*(i^*([f(Y)\times g(Y')])\cdot\s)=i_*([f(Y)]\cdot[g(Y')]\cdot\s)$ and the result follows.
\end{proo}

\begin{coro}
    Let $Y$ be a closed irreducible subvariety of $\G_Q(p,V)$
  (resp. $\G_\omega(p,V)$) such that $[Y]\cdot [Y]\cdot \s\neq0$, then
  $Y$ is algebraically simply connected.
\end{coro}

\begin{proo}
  Follow the proofs of Theorem 5.1 and Corollary 5.3 in \cite{fulton-lazarsfeld}.
\end{proo}

Let us state the following corollary that we shall need in the next section.

\begin{coro}
\label{coro-utile}
  Let $Y$ be a cumbersome irreducible subvariety of $\G_Q(p,V)$ (resp. $\G_\omega(p,V)$) with $\dim Y-\codim Y\geq p$.

(\i) Then for $g$ general in $G$ the intersection $Y\cap (g\cdot Y)$ is irreducible.

(\i\i) Then $Y$ is algebraically simply connected.

(\i\i\i) Then for $X^P(w)$ a Schubert variety of codimension $d$, the intersection $Y\cap (g\cdot X^P(w))$ is irreducible as soon as $d\leq \dim Y-p$ and $\s(w)\cdot\s\neq0$.
\end{coro}

\begin{proo}
  For \textit{(\i)} and \textit{(\i\i)}, we only need to remark that we can apply the former two corollaries because by assumption we have $[Y]\cdot [Y]\cdot \s\neq0$. For \textit{(\i\i\i)}, we need the condition $[Y]\cdot \s(w)\cdot \s\neq0$ which will be satisfied as soon as $\s(w)\cdot\s\neq0$. 
\end{proo}

\section{Transplanting the Picard group}

In this section, we prove transplanting theorems from $X$ to a smooth subvariety $Y$ in $X$ for the groups $N^1(X)$ and $N^1(Y)$ of divisors modulo numerical equivalence. 

\subsection{First method}

Recall the definition of a cumbersome subvariety (Definition \ref{defi-cumbersome}) and remark that this definition only depens on the cohomology class of the subvariety. We will therefore also use the word cumbersome for cohomology classes.

\begin{defi}
  Let $\xi$ be a cohomology class in $H^{2k}(X,\Z)$, we define the bilinear form $(\ ,\ )_\xi$ on $H^{\dim X-k}(X,\Z)$ by $(\s,\s')_\xi=\s\cup\s'\cup\xi$.
\end{defi}

Recall that the classes $(\s^P(w))_{w\in W^P}$ of the Schubert
varieties $(X^P(w))_{w\in W^P}$ form a basis of the cohomology
$H^*(X,\Z)$. We define the subset $W^P_d$ of $W^P$ by $W^P_d=\{w\in
W^P\ /\ \deg(\s(w))=d\}$. We will denote by $H$ a general hyperplane
section in $X$ and by $h$ its cohomology class.

\begin{theo}
\label{theo-transplanting}
  Let $Y$ be a smooth subvariety of $X$ and assume that the following conditions hold:

(H1) there exists $\Xi\subset X$ such that if $\xi=[\Xi]$, then $s=\dim Y-\codim Y-\deg\xi-2$ is even, non negative and the bilinear form $(\ ,\ )_\xi$ is positive definite;

(H2) for $g$ and $g'$ general in $G$, the intersection $Y\cap (g\cdot Y)\cap (g'\cdot \Xi)$ is irreducible;

(H3) there exists an integer $d\in[0,\dim Y-\codim Y-2]$ such that the class $h^d\cdot [Y]$ is cumbersome;

(H4) for all $w\in W^P_{d+\codim Y}$ and for $g$ general in $G$, the intersection $Y\cap (g\cdot X^P(w))$ is irreducible;

Then $N^1(Y)=\Z$.
\end{theo}

\begin{proo}
Let $H$ be a general hyperplane section in $X$ and let us denote by
$H_Y$ its intersection with $Y$. Let us first remark that, replacing
$D$ with $D+m H_Y$ with $m$ large enough, we may and will assume that
$D$ is smooth. We shall denote by $i$ the embedding of $Y$ in $X$ and
set $r=\dim Y-\codim Y-2-d$.

We want to compare $[D]$ and $[H_Y]=i^*h$ with respect to numerical
equivalence. For this, consider in $Y$ the classes of curves 
$c(u)=i^*(\sigma(u)h^{r})$ and $\gamma(v)=[D]i^*(\sigma(v)h^{r})$ for
$u\in W^P_{d+\codim Y+1}$ and $v\in W^P_{d+\codim Y}$. We define the
matrix $M$ with two lines and columns indexed by the disjoint union $W^P_{d+\codim
  Y+1}\cup W^P_{d+\codim Y}$ by 
$$\displaystyle{M= \left(
\begin{array}{cc}
i^*(h)c(u)_{u\in W^P_{d+\codim Y+1}}& i^*(h)\gamma(v)_{v\in W^P_{d+\codim Y}}\\ 

[D] c(u)_{u\in W^P_{d+\codim Y+1}}& [D]\gamma(v)_{v\in W^P_{d+\codim Y}}\\
\end{array}\right)}.$$
Let us note that Lefschetz's hyperplane Theorem is valid for the group
of divisors modulo numerical equivalence (this comes from the fact
that for divisors, being numerically trivial is the same as being of
torsion in the cohomology, see \cite[Section 19.3.1]{fulton}). This
remark together with the same proof as in \cite{arrondo-caravantes}
lead to the following lemma.

\begin{lemm}
\label{lemm-rang1}
  The divisor $D$ and $H_Y$ are colinear in $N^1(Y)$ if and only if $M$ is of rank one.
\end{lemm}

\begin{proo}
  If $D$ is numerically equivalent to a multiple of $H_Y$, then the
  rank of the matrix $M$ is one. Conversely, if the second line of the
  matrix is $q$ times the first one, let $D_0=[D]-qi^*(h)$. The
  intersections $D_0c(u)$ and $D_0\gamma(v)$ for $u\in W^P_{d+\codim
    Y+1}$ and $v\in W^P_{d+\codim Y}$ vanish. But $h^{d+\codim Y+1}$
  and $h^{d+\codim Y}$ are respectively linear combinations of $\s(u)$ and $\s(v)$
  respectively thus $h^{\dim Y-1}$ and $h^{\dim Y-2}$ are respectively
  linear combinations of $\s(u)h^r$ and $\s(v)h^r$ with $u\in
  W^P_{d+\codim Y+1}$ and $v\in W^P_{d+\codim Y}$. Hence we have
  $[D_0]i^*(h)^{\dim Y-1}=0$ and $[D_0][D]i^*(h)^{\dim Y-2}=0$ If $S$ is a
  smooth surface obtained from $Y$ by $\dim Y-2$ hyperplane sections
  and if $j:S\to Y$ is the embedding, we get the equalities
  $j_*(j^*([D_0]i^*(h)))=[D_0]i^*(h)\cap j_*[S]=[D_0]i^*(h)^{\dim
    Y-1}=0$ and $j_*j^*[D_0]^2=[D_0]([D]-q^*(h)))i^*(h)^{\dim
    Y-2}=0$. Thus, in the surface $S$, we have the intersection equalities
$$j^*[D_0]j^*i^*(h)=0 \textrm{ and } j^*[D_0]^2=0.$$
By Hogde index Theorem, the class $j^*[D_0]$ has to be numerically
trivial. By Lefschetz's hyperplane Theorem, this has to be true for
$[D_0]$.
\end{proo}

We are therefore left to prove that the matrix $M$ has rank one. We
first prove that the ``left part'' of $M$, \emph{i.e.} the submatrix
formed by the columns of $M$ indexed by $W^P_{d+\codim Y+1}$, is of
rank one. For this, remark that it is enough to prove that $i_*[D]$
and $i_*i^*(h)=[Y]h$ are colinear in $H^{\codim Y+1}(X,\Z)$. Indeed, the values of the
``left part'' of $M$ can be computed using intersection in
$X$. Namely, we have the equalities $i_*(i^*(h)c(u))=i_*i^*(\s(u)h^{r+1}) =
[Y]\s(u)h^{r+1}$ and $i_*([D]c(u))=i_*([D]i^*(\s(u)h^r)=i_*[D]\s(u)h^r$.

\begin{lemm}
  The classes $i_*[D]$ and $h[Y]$ are colinear in $H^{\codim Y+1}(X,\Z)$.
\end{lemm}

\begin{proo}
Because the subvarieties $Y$ and $D$ are smooth, we have an exact
sequence of normal bundles $0\to N_{D/Y}\to N_{D/X}\to
N_{X/Y}\vert_D\to0$. Taking the top Chern classes and denoting by
$j:D\to Y$ the inclusion, we get the equality 
\begin{equation}
  \label{equa-debase}
  j^*i^*i_*[D]=j^*[D]\cdot j^*(i^*[Y]).
\end{equation}
Let $S_Y$ be the surface obtained from $Y\cap (g\cdot Y)\cap (g'\cdot
\Xi)$ with $g$ and $g'$ general in $G$ by $\dim Y-\codim Y - \deg\xi-2$
general hyperplane sections. This surface $S_Y$ is irreducible by
assumption. Denote by $j_Y:S_Y\to Y$ the embedding. Pushing forward
the equality (\ref{equa-debase}) by $i_*j_*$ we have 
$$i_*[D]^2=i_*((i^*i_*[D])[D])=i_*j_*j^*i^*i_*[D]=i_*j_*j^*([D]i^*[Y])
=i_*([D]i^*[Y][D])=i_*([D][D]i^*[Y]).$$
Let $s={\dim Y-\codim Y-\deg\xi-2}$, multiplying by $\xi h^{s}$ we obtain the equalities
$$\!\!\!
\begin{array}{ll}
(i_*[D],i_*[D])_{\xi h^{s}}\!\!\!\!
&=i_*[D]i_*[D]\xi h^{s}=
i_*([D][D]i^*[Y])\xi h^{s}=
i_*([D][D]i^*([Y]\xi h^{s}))= 
i_*([D][D]{j_Y}_*[S_Y])\\ 
(i_*[D],i_*[D])_{\xi h^s}\!\!\!\!
&=i_*{j_Y}_*(j_Y^*[D]j_Y^*[D])= i_*{j_Y}_*(j_Y^*[D]^2).\\
\end{array}$$
On the one hand, the surface $S_Y$ being irreducible, we may apply
Hodge index theorem on $S_Y$ to get the inequality $j_Y^*[D]^2\cdot
(j_Y^*i^*h)^2\leq (j_Y^*[D]j_Y^*i^*h)^2$. On the other hand, we may
compute the equalities 
$$
\begin{array}{ll}
i_*{j_Y}_*((j_Y^*i^*h)^2)=i_*((i^*h)^2[S_Y])=i_*((i^*h)^2i^*([Y]\xi h^s))=
h^2[Y]\xi [Y]h^s=(h[Y],h[Y])_{\xi h^s} \ \ \textrm{and}\\
i_*{j_Y}_*(j_Y^*[D]j_Y^*i^*h)  =i_*([D]i^*h[S_Y])= i_*(i^*h[D]i^*([Y]\xi h^s))=
i_*[D]h[Y]\xi h^s=(i_*[D],h[Y])_{\xi h^s}.\\
\end{array}$$
All together we get the inequality
$(i_*[D],i_*[D])_{\xi h^s}(h[Y],h[Y])_{\xi h^s}\leq(i_*[D],h[Y])_{\xi h^s}$. But $(\ ,\
)_\xi$ being positive definite, the same is true for $(\ ,\ )_{\xi h^s}$ and by Cauchy-Schwartz we must have equality and the fact that $i_*[D]$ and $h[Y]$ are colinear.
\end{proo}

We are left to prove that the ``right part'' of the matrix $M$,
\emph{i.e.} the submatrix formed by the columns of $M$ indexed by
$W^P_{d+\codim Y}$ is spanned by the ``left part''. For this, let us
write $[Y]h^d=\sum_{w\in W^P_{d+\codim Y}}a_w\s(w)$. By assumption
$[Y]h^d$ is cumbersome thus we have $a_w>0$ for all $w$. Pushing forward
equation (\ref{equa-debase}) by $j_*$ and multiplying by $i^*h^{d+r}$, we get 
\begin{equation}
  \label{equa-decomp}
  [D]i^*(i_*[D]h^{d+r})=[D]^2i^*([Y]h^{d+r})=
[D]^2\sum_{w\in W^P_{d+\codim Y}}a_wi^*(\s(w)h^{r}).
\end{equation}
Let us denote by $S_w$ the surface obtained by $r$ hyperplane intersections of $Y\cap (g\cdot X^P(w))$ for a general $g$ in $G$ and $w\in W^P_{d+\codim Y}$. By assumption $S_w$ is irreducible. We denote by $j_w:S_w\to Y$ the inclusion. We have the equality ${j_w}_*j_w^*([D]^2)=[D]^2i^*(\s(w)h^{r})$ and by Hodge index Theorem on $S_w$, we get $j_w^*([D]^2)\cdot j_w^*i^*(h^2)\leq (j_w^*[D]j_w^*i^*h)^2$. 
Pushing forward equation (\ref{equa-decomp}) with $i_*$, we get
$$(i_*[D])^2h^{d+r}\leq \sum_{w\in W^P_{d+\codim Y}}a_w \frac{i_*{j_w}_*((j_w^*[D]j_w^*i^*h)^2)}{i_*{j_w}_*(j_w^*i^*(h^2))}.$$
But, applying projection formula, we have the two equalities $i_*{j_w}_*(j_w^*i^*(h^2))= [Y]\s(w)h^{r+2}$ and $i_*{j_w}_*((j_w^*[D]j_w^*i^*h)^2)=(i_*[D]\s(w)h^{r+1})^2$. Furthermore, there exists a $\lambda$ such that $i_*[D]=\lambda h[Y]$. We thus have inequality:
$$\lambda^2 [Y]^2 h^{d+r+2} \leq \sum_w
a_w \frac{(\lambda [Y]\s(w)h^{r+2})^2}{[Y]\s(w)h^{r+2}}=\sum_w
a_w \lambda^2 [Y]\s(w)h^{r+2}=\lambda^2 [Y]^2 h^{d+r+2}.$$
We thus have equality in all the above inequalities. In particular, we have equality in the Hodge index inequality for the surface $S_w$: $j_w^*([D]^2)\cdot j_w^*i^*(h^2)= (j_w^*[D]j_w^*i^*h)^2$. Pushing forward with ${j_w}_*$, we get $([D]^2i^*(\s(w)h^r))\cdot(i^*(\s(w)h^{r+2}))=([D]i^*(\s(w)h^{r+1}))^2$. Let us write $h\s(w)=\sum_u c_{w,h}^u\s(u)$ where the sum runs over $u\in W^P_{d+\codim Y+1}$. We get the equality 
$$([D]\gamma(w))\cdot\left(\sum_u c_{w,h}^u i^*(h) c(u)\right)=(i^*(h)\gamma(w))\cdot\left(\sum_u c_{w,h}^u [D]c(u)\right).$$
This can be rewritten as
$$\sum_u c_{w,h}^u\left\vert
  \begin{array}{cc}
    i^*(h)c(u) & i^*(h)\gamma(w) \\

    [D]c(u) & [D]\gamma(w) \\
  \end{array}\right\vert =0.$$
Because for $w$ fixed, all the Littlewood-Richardson coefficients $c_{w,h}^u$ are non negative and at least one of them is positive, one of these minors must vanish completing the proof of the theorem.
\end{proo}

\begin{rema}
  A very natural choice for the subvariety $\Xi$ and the class $\xi$ is to take a complete intersection in $X$. In that case $\xi$ is a multiple of $h^k$ for some interger $k$. 

\begin{fact}
  The value of ${\rm eff}(X)$ is $\dim X-4$ for all rational homogeneous spaces with Picard rank one except for the one given in the following list. 
  
\vskip 0.2 cm

\begin{tabular}{|c|c||c|c||c|c||c|c|}
    \hline
$X$ & {\rm eff}(X) & $X$ & {\rm eff}(X) & $X$ & {\rm eff}(X) & $X$ & {\rm eff}(X) \\
\hline
 $\p^{2n}$ & $0$ & $\G(2,n)$ & $0$ & $E_6/P_1$ & $0$ & $E_8/P_7$ & $\dim X-8$\\
 $\p^{2n-1}$ & $1$ &  $\G_\omega(2,2n)$ & $1$ & $E_7/P_1$ & $1$ & $E_8/P_8$ & $1$\\
$\Q^{2n-1}$ & $1$ &  $\G_\omega(2,2n)$ & $1$ & $E_7/P_6$ & $\dim X-8$ & $F_4/P_1$ & $1$ \\
$\Q^{4n}$ & $0$ &  $\G_Q(2,4n+2)$ & $4n+1$ & $E_7/P_7$ & $\dim X- 8$ & $F_4/P_4$ & $1$ \\
$\Q^{4n+2}$ & $2$ & $\G_Q(2,4n)$ & $1$ & $E_8/P_1$ & $\dim X-12$ & & \\
\hline
  \end{tabular}
\end{fact}
Remark that we will be interested in the homogeneous spaces $X$ with small values of ${\rm eff}(X)$. Among those we have the projective spaces over composition algebras ($\p^n$, $\G(2,n)$ and $E_6/P_1$) and also adjoint and coadjoint varieties (for example the varieties of isotropic line $\G_Q(2,n)$ and $\G_\omega(2,2n)$). For adjoint and coadjoint varieties, the description of the cohomology classes in terms of roots given in \cite{adjoint} easily implies that if the Weyl involution of the Dynkin diagram is trivial, then ${\rm eff}(X)=1$. This occurs in all types except type $A_n$, type $D_{2n+1}$ and type $E_6$.

In particular, to prove Theorem \ref{theo-main3}, we only need to deal with the varieties $X$ with ${\rm eff}(X)=0$ or ${\rm eff}(X)=1$ because the other cases follow from Sommese and van de Ven result: Theorem \ref{theo-sommese-van-de-ven}.
\end{rema}

We can now state the following corollary of Theorem \ref{theo-transplanting}.

\begin{coro}
\label{coro-classique}
  Let $X$ be $\G_Q(2,2n+1)$, $\G_\omega(2,2n)$ or $\G_Q(2,4n)$ and $Y$ be a smooth cumbersome subvariety with $2\dim Y\geq{\dim X+2}$, then $\pic(Y)=\Z$.
\end{coro}

\begin{proo}
   We have ${\rm eff}(X)=1$ thus we can choose for $\Xi$ an hyperplane
  section and take $d=0$ (note that the dimension of $X$ is odd thus
  we have $2\dim Y\geq{\dim X+3}$). Using the results in
  \cite{adjoint}, any Schubert class corresponds to a root of the
  group $G$. The class $\s$, with notation as in Corollary
  \ref{coro-utile}, corresponds to the root $\Theta-\a_1-\a_2$ where
  $\Theta$ is the highest root (resp. highest short root) if the form
  is a symmetric form $Q$ (resp. a symplectic form $\omega$) and
  $\a_1$ and $\a_2$ are simple roots with notation as in Bourbaki
  \cite{bou}. The Schubert class $\s'$ of minimal degree with $\s\cdot
  \s'=0$ corresponds to the root $-\a_1$ and its degree satisfies
  $2\deg\s'>\dim X$. Thus we have $\s\cdot \s(w)\neq0$ for
  $2\deg\s(w)\leq\dim X$. Therefore we have $[Y]\s\neq0$ and for $w\in
  W^P_{\codim Y}$ we have $\s(w)\s\neq0$. Because $Y$ is cumbersome we
  get $[Y]^2\s\neq0$ and $[Y]\s(w)\s\neq0$. Corollary \ref{coro-utile}
  implies that $Y\cap (g\cdot Y)$ and $Y\cap (g\cdot X^P(w))$ are
  irreducible for $g$ general in $G$ and for $w\in W^P_{\codim Y}$ and
  that $Y$ is algebraically simply connected. By Theorem
  \ref{theo-transplanting}, any divisor in $Y$ is numerically
  equivalent to the hyperplane section and because $Y$ is
  algebraically simply connected, the result follows.
\end{proo}

\begin{rema}
  The same method also proves that for $Y$ smooth cumbersome subvariety of $X$ with $2\dim Y\geq \dim X+2$ (resp. $2\dim Y\geq 2\dim X+4$) where $X$ is a projective space, a smooth quadric with $\dim X\not\equiv 2\ ({\rm mod}\ 4)$ (resp. a smooth quadric with $\dim X\equiv 2\ ({\rm mod}\ 4)$) or a grassmannian $\G(2,n)$, then $\pic(Y)=\Z$. Remark also that the assumption cumbersome can be dropped when $X$ is a projective space or a quadric because in these dimensions all subvarieties in $X$ are cumbersome.

These result are already known. For the projective space this is part of Barth and Larsen \cite{barth-larsen} results. The results for quadrics follow from Barth-Larsen's results (except for $2\dim Y=\dim X+2$). The last case for quadrics and the case of the grassmannian $\G(2,n)$ are proved in \cite{arrondo-caravantes}.
\end{rema}

\begin{rema}
\label{rema-arr-car}
  (\i) The result of Corollary \ref{coro-classique} for $X=\G_\omega(2,2n)$ can be deduced from the corresponding result for $X=\G(2,2n)$ proved by Arrondo and Caravantes in \cite{arrondo-caravantes}. Indeed, the variety $\G_\omega(2,2n)$ is an hyperplane section of $\G(2,2n)$ therefore if $Y$ is a cumbersome subvariety in $\G_\omega(2,2n)$, then it is also cumbersome in $\G(2,2n)$.

(\i\i) For $X=\G_Q(2,n)$ however, the result of Arrondo and Caravantes do not give any information. Indeed, let $Z_v=\{V_2\in\G(2,n)\ /\ V_2 \textrm{ contains a fixed vector $v$}\}$. Then if $v$ is non isotropic, the intersection $Z_v\cap \G_Q(2,n)$ is empty \emph{i.e.} $[Z_v]\cdot[\G_Q(2,n)]=0$. In particular, for any subvariety $Y$ of $\G_Q(2,n)$, we have $[Y]\cdot[Z_v]=0$ thus if $2\dim Y\geq\dim\G(2,n)+2$, then $\dim Y\geq\codim Z_v=\frac{1}{2}\dim\G(2,n)$ and $Y$ is not cumbersome in $\G(2,n)$.
\end{rema}

\subsection{Arrondo-Caravantes technique}

In this section we will use a more direct generalisation of E. Arrondo and J. Caravantes technique to prove the same type of results on the Picard group of smooth subvarieties of small codimension in homogeneous spaces but without the assumption \textit{(H2)} in Theorem \ref{theo-transplanting}. This will be useful for exceptional groups for which the results of Section \ref{section-connexite} do not apply. Next theorem in particular completes the proof of Theorem \ref{theo-main3}. However, the technique here is more \emph{ad-hoc} than in the previous section. In particular it does not clearly explain why the fact that the form $(\ ,\ )_{h^k}$ for some $k$ is positive definite should be important.

\begin{theo}
\label{theo-trans-bis}
  Let $X$ be an adjoint or coadjoint variety with ${\rm eff}(X)=1$. If $Y$ is a smooth cumbersome subvariety in $X$ with $2\dim Y\geq \dim X+2$, then $N^1(Y)=\Z$.
\end{theo}

\begin{rema}
  The adjoint or coadjoint varieties with ${\rm eff}(X)=1$ are $\Q^{2n-1}$, $\G_Q(2,2n+1)$, $\G_Q(2,4n)$, $\p^{2n-1}$, $\G_\omega(2,2n)$, $E_7/P_1$, $E_8/P_8$, $F_4/P_1$, $F_4/P_4$, $G_2/P_1$ and $G_2/P_2$.
\end{rema}

\begin{proo}
  We shall present the proof in general but we shall only prove the computational Lemma \ref{lemm-def-pos} for exceptional types. For classical groups, the result of the theorem was proved in Corollary \ref{coro-classique}. We start as in the proof of Theorem \ref{theo-transplanting} and use the same notation as in that proof. We set $d=1$. By taking general hyperplane sections, we may assume that $\dim Y$ is minimal in the range that is (because $\dim X$ is odd) we have the equality $2\dim Y= \dim X+3$. We want to prove that the rank of the matrix $M$ is one. Keeping notation as in the proof of Theorem \ref{theo-transplanting}, we have $r=0$ and for $w\in W^P_{\codim Y+1}$, we define the surface $S_w$ to be the intersection $Y\cap(g\cdot X^P(w))$ for $g$ general in $G$.

 \begin{lemm}
    The surfaces $S_w$ for $w\in W^P_{\codim Y+1}$ are irreducible.
  \end{lemm}  

\begin{proo}
Using the results in \cite{adjoint}, if $X$ is adjoint (resp. coadjoint), then any Schubert class $\s(w)$ is represented by a long (resp. short) root $\a_w$ and any long (resp. short) root is associated to a Schubert class. It is easy to check that $\Sigma^P(w)=\{\beta \textrm{ simple root}\ /\ \sca{\beta^\vee}{\a_w}>0\}$. For $w\in W^P_{\codim Y+1}$, the root $\a_w$ is simple root. Take $v\in W^P_{\codim Y+2}$ such that $\a_{v}=-\a_w$, then $\Sigma^P(v)\cap\Sigma^P(w)=\emptyset$. Furthermore, the unipotent subgroup $U(\a_w)$ is contained in $S^P(w)$ thus $S^P(w)\cdot X^P(v)$ contains $X^P(s_{\a_w}(v))$ and by the results of \cite{adjoint} the associated root is $\a_{s_{\a_w}(v)}=s_{\a_w}(\a_v)=\a_w$. Thus $S^P(w)\cdot X^P(v)=X^P(w)$. Because $X^P(v)$ is of codimension $\codim Y+2$ and $Y$ is cumbersome, we have $[Y]\cdot \s(v)\neq0$ and we may apply Theorem \ref{theo-bertini} to finish the proof. 
\end{proo}

The same computation as in the proof of Theorem
\ref{theo-transplanting} gives us the inequality
$$q(i_*[D]):=(i_*[D])^2h - \sum_w a_w
\frac{(i_*[D]\s(w)h)^2}{[Y]\s(w)h^{2}}\leq0.$$
The above expression $q(i_*[D])$ is a quadratic form in $i_*[D]$. We
shall compute it explicitely. For this we first remark that the
elements $u\in W^P_{\codim Y+2}$ are in correspondence (using the
results in \cite{adjoint} once again) with the opposite of simple roots
\emph{i.e.}, if $X$ is adjoint (resp. coadjoint), the root $-\a_u$ is a
simple long (resp. short) root. For a pair of opposite of simple roots
$(\a_u,\a_{u'})$ with $\sca{\a_u^\vee}{\a_{u'}}<0$, we define the coordinate
$$L(u,u')=\frac{1}{\sqrt{i^*(h)c(u) \cdot i^*(h)c(u')}}\left\vert
\begin{array}{cc}
i^*(h)c(u) & i^*(h)c(u')\\

[D]c(u) & [D]c(u')\\ 
\end{array}\right\vert.$$
We shall prove in Lemma \ref{lemm-def-pos} that $q(i_*[D])$ can be
writen as a quadratic form ${\cal Q}(L(u,u')_{u,u'})$ in the variables
$L(u,u')$ and that ${\cal Q}(L(u,u')_{u,u'})$ is positive definite. As
we have ${\cal Q}(L(u,u')_{u,u'})\leq0$, this implies that we have
equality and that $L(u,u')=0$ for all $(u,u')$ with
$\sca{\a_u^\vee}{\a_{u'}}<0$. In particular, this implies that the
``left part'' of the matrix $M$ is of rank one and by the same
argument as in the proof of Theorem \ref{theo-transplanting}, that the
``right part'' of $M$ is colinear to its ``left part''. Thus $M$ is of
rank one and the result follows by Lemma \ref{lemm-rang1}.
\end{proo}

To finish the proof of Theorem \ref{theo-trans-bis} we need to compute more explicitely the quadratic form $q(i_*[D])$ defined above. For this we shall use extensively the correspondence between roots and cohomology classes in an adjoint or coadjoint variety $X$. First of all we note that for $Y$ with $2\dim Y=\dim X+3$ and for $w\in W^P_{\codim Y+1}$, the associated root $\a_w$ is simple while for $u\in W^P_{\codim Y+2}$, the associated root $\a_u$ is the opposite of a simple root. Furthermore, if ${\rm eff}(X)=1$, then the Poincar{\'e} duality maps $\s(w)$ with $w\in W^P_{\codim Y+1}$ to $\s(u)$ with $u\in W^P_{\codim Y+2}$ and $\a_u=-\a_w$. To simplify notation, we shall identify $w\in W^P_{\codim Y+1}$ with the simple root $\a_w$ and $u\in W^P_{\codim Y+2}$ with the opposite simple root $\a_u$. Therefore it will make sense to write $-w$ or $-u$, for example $\s(-w)=\s(u)$ for $\a_u=-\a_w$. We also set $x_u=i^*(h)c(u)$, $y_u=[D]c(u)$ and $d_u=[Y]\s(-u)h^2$ for $u\in W^P_{\codim Y+2}$. With these notation, we have $x_u=i^*(h\s(u))=h[Y]\s(u)=a_{-u}$.

\begin{lemm}
  \label{lemm-def-pos}
We assume that the group $G$ is exceptional. Let us choose a complete order on the set of simple roots (for example the one given in \cite{bou}). Let $\cal P$ be the set of couple of simple roots $(\a,\beta)$ with $\a<\beta$ and $\sca{\a^\vee}{\beta}<0$. We define the matrix ${\cal Q}=(q_{(u_1,u'_1),(u_2,u'_2)})_{((u_1,u'_1),(u_2,u'_2))\in{\cal P}^2}$ with index set $\cal P$ by
$$q_{(u_1,u'_1),(u_2,u'_2)}=\left\{
  \begin{array}{ll}
1-\frac{x_{u_1}}{d_{u_1'}}-\frac{x_{u'_1}}{d_{u_1}} & \textrm{for $(u_1,u'_1)=(u_2,u'_2)$},\\
\frac{\sqrt{x_{u_1}x_{u'_2}}}{d_{u_1'}} & \textrm{for $u'_1=u_2$},\\
0 & \textrm{otherwise.}\\
  \end{array}
\right.$$

(\i) Then the quadratic form ${\cal Q}$ is positive definite.

(\i\i) We have the formula $q(i_*[D])={\cal Q}(L(u,u'))_{(u,u')\in{\cal P}}$.
\end{lemm}

\begin{proo}
  \textit{(\i)} The first statement is only computational once we remark that $d_u$ can be expressed with the $x_u$: we have $d_u=h[Y]\s(-u)h=\sum_{u'}c_{-u,h}^{u'}h[Y]\s(u')=\sum_{u'}c_{-u,h}^{u'}x_{u'}$. We therefore need to check that in all cases, the above quadratic form is positive definite and we can check this by computing the principal minors. An easy computation gives the result.

Remark that for classical case, the computation is less easy since we get matrices indexed by $n$. For example, in type $B_n$ (resp. $C_n$) for long (resp. short) roots, an easy induction gives that the determinant of the matrix defining ${\cal Q}$ is
$$\frac{nx_{u_1}\cdots x_{u_{n-1}}}{d_{u_1}\cdots d_{u_{n-1}}}$$
where $(u_i)_{i\in[1,n-1]}$ correspond to the opposite of the simple long (resp. short) roots $(\a_i)_{i\in[1,n-1]}$ with notation as in \cite{bou}. For type $D_{2n}$, we do not know such a simple formula but we expect that the quadratic form ${\cal Q}$ is also positive definite. By a direct check we proved this for type $D_4$.

  \textit{(\i\i)} Let us define the following quadratic form:
$$q'=\sum_{u<u'}\sum_{u''}\frac{x_{u''}}{d_{u''}}c_{-u'',h}^uc_{-u'',h}^{u'}L(u,u')^2 -\sum_{u<u'}c_{-u,h}^{u'}L(u,u')^2.$$

\begin{fact}
  We have the equality $q=q'$.
\end{fact}

\begin{proo}
  We compute the coefficient of $y_uy_{u'}$ for $u$ and $u'$ in $W^P_{\codim Y+2}$. We have $i_*[D]=\sum_uy_u\s(-u)$ therefore we get $(i_*[D])^2h=\sum_{u,u'}c_{-u,h}^{u'}y_uy_{u'}$. We also have for $w=-u''$ the equality $i^*[D]\s(w)h=\sum_{u}c_{-u'',h}^{u}y_{u}$. We thus get
$$q(i_*[D])=\sum_{u,u'}c_{-u,h}^{u'}y_uy_{u'}-\sum_{u,u'}\left(\sum_{u''}\frac{x_{u''}}{d_{u''}}c_{-u'',h}^{u}c_{-u'',h}^{u'}\right)y_uy_{u'}.$$
where $u$, $u'$ and $u''$ run in $W^P_{\codim Y+2}$. Let us set
$$A_{u,u'}=c_{-u,h}^{u'}-\sum_{u''}\frac{x_{u''}}{d_{u''}}c_{-u'',h}^{u}c_{-u'',h}^{u'} \textrm{ so that } q=\sum_{u,u'}A_{u,u'}y_uy_{u'} \textrm{ and } 
q'=-\sum_{(u,u')\in{\cal P}}A_{u,u'}L(u,u')^2.$$
Let us note that $L(u,u')^2=\frac{x_{u'}}{x_{u}}y_u^2+\frac{x_{u}}{x_{u'}}y_{u'}^2-2y_uy_{u'}$. The coefficient $q'_{u,u'}$ of $y_uy_{u'}$ in $q$ is therefore given by 
$$q'_{u,u'}=\left\{
  \begin{array}{ll}
\displaystyle{-\sum_{u''\neq u}\frac{x_{u''}}{x_u}A_{u,u''}}& \textrm{if $u=u'$}\\
\\
\displaystyle{2A_{u,u'}}&\textrm{if $u\neq u'$}\\
  \end{array}
\right.$$
Thus if $q_{u,u'}$ is the coefficient of $y_uy_{u'}$ in $q$, we clearly have $q_{u,u'}=q'_{u,u'}$ for $u\neq u'$. But we can also compute, using the formulas $d_{u'}=\sum_{u''}c_{-u',h}^{u''}x_{u''}$ and $c_{-u',h}^{u}=c_{-u,h}^{u'}$, the equalities:
$$
\begin{array}{ll}
q'_{u,u}&=\displaystyle{-\sum_{u''\neq u}\frac{x_{u''}}{x_u}c_{-u,h}^{u''}+ \sum_{u'}\frac{x_{u'}}{d_{u'}x_u}c_{-u',h}^{u}\sum_{u''\neq u}c_{-u',h}^{u''}x_{u''}}\\
\\
& =\displaystyle{ -\sum_{u''\neq u}\frac{x_{u''}}{x_u}c_{-u,h}^{u''}+\sum_{u'}\frac{x_{u'}}{d_{u'}x_u}c_{-u',h}^{u}(d_{u'}-c_{-u',h}^hx_u)}\\
\\
& =\displaystyle{ -\sum_{u''\neq u}\frac{x_{u''}}{x_u}c_{-u,h}^{u''}+\sum_{u'}\frac{x_{u'}}{x_u}c_{-u',h}^{u}- \sum_{u'}\frac{x_{u'}}{d_{u'}x_u}c_{-u',h}^{u}c_{-u',h}^hx_u}\\
\\
& =\displaystyle{ c_{-u,h}^u- \sum_{u'}\frac{x_{u'}}{d_{u'}}(c_{-u',h}^{u})^2=A_{u,u}.}\\
\end{array}
$$
This completes the proof of the fact.
\end{proo}

Now we only need to compute the matrix of quadratic form $q'$ which is already in the variable $L(u,u')$. But note that we only want to deal with pairs of simple roots $(\a,\beta)$ or associated elements $(u,u')$ in $W^P_{\codim Y+2}$ with $u<u'$ and $\sca{\a^\vee}{\beta}<0$ or with our notation $\sca{u^\vee}{u'}<0$. We first recall the possible values for $c_{-u,h}^{u'}$:
$$c_{-u,h}^{u'}=\left\{
  \begin{array}{ll}
2 & \textrm{if $u=u'$,}\\
1 & \textrm{ if $\sca{u^\vee}{u'}<0$,}\\
0 & \textrm{ otherwise.}\\
  \end{array}\right.$$
In particular, we see that the only factors $L(u,u')$ appearing in $q'$ with $\sca{u^\vee}{u'}\geq0$ are such that there exists $u''$ with $\sca{u^\vee}{u''}<0$ and $\sca{{u''}^\vee}{u'}<0$. The coefficient of $L(u,u')^2$ being in that case $x_{u''}/d_{u''}$. But because of the obvious identity
$$\left\vert
\begin{array}{ccc}
  x_u & x_{u'} & x_{u''} \\
  y_u & y_{u'} & y_{u''} \\
  x_u & x_{u'} & x_{u''} \\
\end{array}\right\vert=0$$
we have the formula $\sqrt{x_{u''}}L(u,u')=\sqrt{x_{u'}}L(u,u'')+\sqrt{x_{u}}L(u'',u')$. Replacing the factors $L(u,u')$ with $\sca{u^\vee}{u'}\geq0$ using this formula, we can write $q'$ only with factors $L(u,u')$ with $(u,u')\in {\cal P}$. We are then left to compute the matrix of this quadratic form. Now the coefficient in $q'$ of a factor $L(u,u')^2$ with $(u,u')\in{\cal P}$ has for contribution $-c_{-u,h}^{u'}=-1$ from the right hand side of $q'$, has, from the left hand side when $u''=u$ or $u''=u'$, the contribution $2x_u/d_u+2x_{u'}/d_{u'}$ and its last contribution comes from the left hand side for $u''\in {\cal B}(u,u')=\{u''\neq u'\ /\ \sca{{u''}^\vee}{u}<0\}$ or $u''\in {\cal B}(u',u)=\{u''\neq u\ /\ \sca{{u'}^\vee}{u''}<0\}$ and is equal to
$$\sum_{u''\in{\cal B}(u,u')}\frac{x_{u}}{d_u}\frac{x_{u''}}{x_u}+\sum_{u''\in{\cal B}(u',u)}\frac{x_{u'}}{d_u'}\frac{x_{u''}}{x_u'}.$$
Summing these contributions, using the values of $c_{-u,h}^{u'}$ and the equality $d_u=\sum_{u'}c_{-u,h}^{u'}x_{u'}$, we get the diagonal terms in the matrix ${\cal Q}$. The only non diagonal terms come from the left hand side of $q'$ for pairs $(u,u')$ with $\sca{u^\vee}{u'}\geq0$ and elements $u''$ with $\sca{u^\vee}{u''}<0$ and $\sca{{u''}^\vee}{u'}<0$. We easily get the non diagonal terms of ${\cal Q}$ this way.
\end{proo}

\begin{coro}
  Let $X$ be an adjoint or a coadjoint varieties with ${\rm eff}(X)=1$, then if $Y$ is a smooth cumbersome and simply connected subvariety in $X$ with $2\dim Y\geq \dim X+2$, then $\pic(Y)=\Z$.
\end{coro}

We finish this section with few examples of smooth subvarieties with picard number greater than 1 in rational homogeneous spaces with Picard number one.

\begin{exam}
  Embed $\pu\times\p^{n-1}$ in $\p^{2n-1}$ via the Segre embedding, we get a smooth subvariety of $\p^{2n-1}$ with dimension $n$. If furthermore $n$ is even, then as already noticed in \cite{arrondo-caravantes}, the image of $\pu\times\p^{n-1}$ is contained in a smooth quadric giving an example of a smooth variety $Y$ with $\pic(Y)=\Z^2$, $\dim Y=n$ in a smooth quadric of dimension $2n-2$ for $n$ even.
\end{exam}

\begin{exam}
  Let $V$ be a vector space of dimension $n$, eventually endowed with a non degenerate symmetric or symplectic form $Q$ or $\omega$. Consider a decomposition $V=U\oplus W$ and look at the subvariety $Y$ of the grassmannian $\G(2,V)$ (and its intersection with $\G_Q(2,V)$ or $\G_\omega(2,V)$) defined by
$$Y=\{V_2\in\G(2,V)\ /\ \dim(V_2\cap U)=\dim(V_2\cap W)=1\}.$$
Then $Y$ is a smooth variety isomorphic to $\p^{\dim U-1}\times \p^{\dim W-1}$ and thus of dimension $n-2$ and Picard number 2. Its intersection with $\G_Q(2,V)$ and $\G_\omega(2,V)$ depends on the restriction of the forms $Q$ and $\omega$ on $U$ and $W$. If both restriction are non degenerate and $U$ and $W$ are orthogonal, we get that $Y\cap \G_Q(2,V)$ is a product of two quadrics of dimensions $\dim U-2$ and $\dim W-2$ thus $\dim(Y\cap \G_Q(2,V))=n-4$ and $Y\cap \G_\omega(2,V)=Y$. If on the contrary $U$ and $W$ are isotropic subspaces, then $Y$ is the incidence variety in $\p^{[n/2]-1}\times\p^{[n/2]-1}$ and is of dimension $2[n/2]-3$. We produce in this way examples of smooth subvarieties $Y$ with Picard number 2 in the grassmannians of lines with maximal dimensions given in the following array.
$$
\begin{array}{cccccc}
  X & \dim X & \dim Y \\
\hline 
\G(2,n) & 2n-4 & n-2 \\
\G_Q(2,2n+1) & 4n-5 & 2n-3 \\
\G_\omega(2,2n) & 4n-5 & 2n-2 \\
\G_Q(2,2n) & 4n-7 & 2n-3 \\
\end{array}$$
Some of the above examples can be generalised for adjoint (resp. coadjoint) varieties of the group $G$. We can embed the adjoint variety corresponding to the maximal subgroup of type $A$ in $G$ obtained from the subsystem of long (resp. short) roots. We produce in this way examples of smooth subvarieties $Y$ with Picard number 2.
However, the bounds we obtain this way are far below the bound given in Theorem \ref{theo-trans-bis}.
\end{exam}

\begin{exam}
  There is one example due to the exceptional isomorphism between type $A_3$ and type $D_3$ giving a variety above the bounds in Theorems \ref{theo-transplanting} and \ref{theo-trans-bis} showing that the condition \emph{cumbersome} cannot be removed to easily. Indeed, consider the inclusions $\G_Q(2,6)\subset\G(2,6)$ and $\G_Q(2,6)\subset\G_Q(2,7)$. The variety $\G_Q(2,6)$ is the incidence variety in $\p^3\times{\p^3}^\vee$ thus is smooth with Picard number 2 and dimension 5. The dimensions of $\G(2,6)$ and $\G_Q(2,7)$ are 8 and 7 respectively therefore we are above the bounds in Theorems \ref{theo-transplanting} and \ref{theo-trans-bis}. However, $\G_Q(2,6)$ is neither cumbersome in $\G(2,6)$ nor in $\G_Q(2,7)$
\end{exam}

\section{Results on the topological fundamental group}

In this section we prove that the results of Section \ref{section-connexite} are still valid for non proper morphisms and therefore deduce results for the topological fundamental group. Let us fix some notation: we keep $V$, $Q$, $\omega$ and as in Section \ref{section-connexite}. We take $X=\G_Q(p,V)$ or $X=\G_\omega(p,V)$. Let $W_q$ be a fixed isotropic $q$-dimensional subspace of $V$ and let $W_{q+1}$ be a dimension $q+1$ isotropic subspace of $V$. Denote by $X^P(w)$ the Schubert variety $\{V_p\ /\ V_p\subset W_q^\perp\}$ and by $\s_{q+1}$ the cohomolgy class of the Schubert variety $\{V_p\ /\ V_p\subset W_{q+1}^\perp\}$ of $X$.

To extend the results of Section 2 as annonced, we only need to prove the following proposition.

\begin{prop}
  Let $Y\to X$ be a morphism with $Y$ irreducible and such that $[f(Y)]\s_{q+1}\neq0$.

(\i) The inverse image $f^{-1}(g\cdot X^P(w))$ is irreducible for $g$ in a dense open subset of $G$.

(\i\i) If $Y$ is unibranch, then $\pi_1(f^{-1}(g\cdot X^P(w)))\to\pi_1(Y)$ is surjective for $g$ general in $G$.
\end{prop}

To prove this proposition, we proceed as for Theorem \ref{theo-connexite} and we therefore only need the statement corresponding to Proposition \ref{prop-bertini-min} for non proper maps. This is done in the next lemma and finishes the proof of Theorem \ref{theo-main2}.

\begin{lemm}
   Let $Y\to X$ be a dominant map with $Y$ irreducible and let $L$ be a general line in $X$.

(\i) The inverse image $f^{-1}(L)$ is irreducible.

(\i\i) If $Y$ is unibranch, then $\pi_1(f^{-1}(L))\to\pi_1(Y)$ is surjective.
\end{lemm}

\begin{proo}
Let $r$ is the maximal dimension of isotropic subspaces in $V$. We start to prove that if \textit{(\i)} holds for $p=r$, then it holds for all $p$. We then prove \textit{(\i)} for $p=r$.

Let us consider the flag variety $\F(p,r,V)$ of partial flags
$V_q\subset V_r$ where $V_q$ and $V_r$ are isotropic subspaces of
dimension $p$ and $r$ respectively in $V$. We denote by $p_1$,
resp. $p_2$ the map $(V_p,V_r)\mapsto V_p$, resp. $(V_p,V_r)\mapsto
V_r$. Let $Z$ be the fiber product obtained from $f$ and $p_1$. We
have a dominant map $g:Z\to\F(p,r,V)$ and by composition $h:=p_2\circ
g$ is dominant. Applying our assumption that \textit{(\i)} holds for
$p=r$ we get that the inverse image of a general line by $h$ is
irreducible. But such a line has the form $\{V_r\ /\ V_r\supset
W_{r'}\}$ where $W_{r'}$ is an isotropic subspace of dimension $r'$ in
$V$ with $r'=r-1$ in all cases except for $V$ endowed with a symmetric
bilinear form $Q$ and $\dim V=2r$ where $r'=r-2$. We thus get that for $W_{r'}$ general in $W$, the subvariety 
$$\{(y,V_r)\ /\ f(y)\subset V_r\supset W_{r'}\}$$
is irreducible. We may therefore assume that there exist such a $W_{r'}$ and $y_0\in
Y^{\rm sm}$ with $f(y_0)\subset W_{r'}$. In that case, the map $(y,V_r)\mapsto V_r$ has a section given by $V_r\mapsto (y_0,V_r)$ and by Lemma \ref{lemm-section} we get that the variety $\{y\ /\ f(y)\subset V_r\}$ is irreducible for $V_r$ general. We are then restricted to the case of a grassmannian variety and may apply Debarre's result (see the proof of \cite[Th\'eor\`eme 6.1]{debarre}) to conclude.

We are therefore left to prove \textit{(\i)} for $p=r$. We proceed by
induction on $r$. For $r=1$ the result is clear because $X$ is itself
a line. To prove the induction step, let us fix an isotropic vector
$v$ in $V$ and consider the map $p_v$, defined on an open subset $U_v$
of $X$, by $V_r\mapsto\pi_v(V_r\cap v^\perp)$ where $\pi_v$ is the
projection from $v$. The open subset $U_v$ has a complementary of
codimension at least two therefore a general line is contained in
$U_v$. The image of $p_v$ is the variety $X_v$ of isotropic subspaces
of (maximal) dimension $r-1$ in $v^\perp/v$. We proved in
\cite[Proposition 5]{fourier} that $p_v$ is a sequence of affine
bundles. In particular, the composition $p_v\circ f$ is dominant and by induction,
the inverse image of a general line is irreducible. We therefore have, by
the description of lines in $X_v$ given above, that for a general
isotropic subspace $W_{r'}$ of dimension $r'$ containing $v$ (and thus
contained in $v^\perp$), the variety 
$$Y_{W_{r'}}=\{y\in Y\ /\ \dim(f(y)\cap W_{r'})\geq r'-1\}$$
is irreducible. Letting $v$ vary, this irreducibility is true for a
general isotropic subspace $W_{r'}$. We may therefore assume that
there is $y_0\in Y^{\rm sm}$ and $W_{r'}$ such that $Y_{W_{r'}}$ is
irreducible and $f(y_0)\supset W_{r'}$. The variety $Y_{W_{r'}}$ is
the inverse image of the Schubert variety $X^P(w)=\{V_r\ /\
\dim(V_r\cap W_{r'})\geq r'-1\}$. The singular locus of $X^P(w)$ is
the Schubert variety $X^P(u)=\{V_r\ /\ V_r\supset W_{r'}\}$ which is a
line in $X$. Let us notice that 
general lines in $X^P(w)$ do not meet $X^P(u)$. Indeed, let $V_{r'}$
an isotropic subspace of dimension $r'$ in $V$ such that
$\dim(V_{r'}\cap W_{r'})=r'-1$ and $V_{r'}\not\subset
W_{r'}^\perp$. Then the line $\{V_r\ /\ V_r\supset V_{r'}\}$ is
contained in $X^P(w)$ but does not meet $X^P(u)$. We can therefore
resctrict ourselves to the open subset $X^P(w)\setminus X^P(u)$ which
is isomorphic to an open subset of $I=\{(V_r,V_{r'-1})\ /\ V_r\supset
V_{r'-1}\subset W_{r'}\}$. We have a natural map $\phi$ on $I$
defined by
$V_r\mapsto V_r\cap W_{r'}$. Its image is the projective space
$\p^{r'-1}$ of $r'-1$ subspaces in $W_{r'}$ and the composition $\phi\circ f:
Y_{W_{r'}}\to \p^{r'-1}$ is dominant. Furthermore, the map
$V_{r'-1}\mapsto (y_0,V_{r'-1})$ is a section of $\phi\circ f$. By
Lemma \ref{lemm-section}, the fibre $F$ of this map is irreducible. But
the image of $F$ under $f$ is either $\p^3$ or a smooth quadric of
dimension 3. For both of them we already know that the inverse image
of a line is irreducible (for the quadric, use the first part of the
proof as $\G_Q(2,5)$ is isomorphic to $\p^3$). This finishes the proof
of \textit{(\i)}.

  We obtain \textit{(\i\i)} by applying \textit{(\i)} to the composition $f\circ \pi$ where $\pi:\widetilde{Y}\to Y$ is the universal covering of $Y$.
\end{proo}

\bigskip\noindent
Nicolas {\sc Perrin}, \\
{\it Hausdorff Center for Mathematics,}
Universit{\"a}t Bonn, Villa Maria, Endenicher
Allee 62, 
53115 Bonn, Germany and \\
{\it Institut de Math{\'e}matiques de Jussieu,} 
Universit{\'e} Pierre et Marie Curie, Case 247, 4 place
Jussieu, 75252 Paris Cedex 05, France.

\noindent {\it email}: \texttt{nicolas.perrin@hcm.uni-bonn.de}.

\end{document}